\documentclass[12pt,reqno]{article}
\usepackage{amssymb}
\usepackage{amsmath}

\newcommand{\nick}[1]{{\bf [~Nick June '07:\ } {\em #1}{\bf~]}}
\newcommand{\carlosa}[1]{{\bf [~Carlos June '07:\ } {\em #1}{\bf~]}}
\def\smallpage{
\addtolength\textwidth{2cm}
\addtolength\oddsidemargin{-1cm}
\addtolength\textheight{3cm}
\addtolength\topmargin{-1.5cm}}
\smallpage

 \newcommand{\bean}{\begin{eqnarray*}}
\newcommand{\eean}{\end{eqnarray*}}
 \def\proof{\noindent{\bf Proof}\ \ }
\def\qed{~~\vrule height8pt width4pt depth0pt}
\def\pix{p}

\begin{document}

\newtheorem{thm}{Theorem}[section]
\newtheorem{cor}{Corollary}[section]
\newtheorem{lemma}{Lemma}[section]
\newtheorem{prop}{Proposition}[section]
\newtheorem{alg}{Algorithm}[section]
\newtheorem{defn}{Definition}[section]
\newtheorem{rem}{Remark}[section]
 \numberwithin{equation}{section}

\def\pr{{\bf P}}

\title{Induced forests in regular graphs with large girth}%
\author{
\begin{centering}
\begin{minipage}{5cm}
Carlos Hoppen
 \\
{\small {\tt  choppen@math.uwaterloo.ca}}\end{minipage}
and\qquad
\begin{minipage}{4.5cm}
Nicholas  Wormald\thanks{Research supported by the  Canada Research
Chairs Program
and NSERC.} \\
{\small {\tt  nwormald@uwaterloo.ca }}
\end{minipage}
\end{centering}\\\\
{
Department of Combinatorics and Optimization }
\\  { University of Waterloo}
\\ {Waterloo ON}\\
{Canada N2L 3G1}
}
\date{}
\maketitle

\begin{abstract}
An induced forest of a graph $G$ is an acyclic induced subgraph of
$G$. The present paper is devoted to the analysis of a simple
randomised algorithm that grows an induced forest in a regular
graph. The expected size of the forest it outputs provides a lower
bound on the maximum number of vertices in an induced forest of
$G$. When the girth is large and the degree is at least 4, our
bound coincides with the best bound known to hold asymptotically
almost surely for random regular graphs. This results in an
alternative proof for the random case.
\end{abstract}

\section{Introduction}

An \textit{induced forest} in a graph $G$ is an acyclic induced
subgraph of $G$. The problem of finding a large induced forest in
a graph $G$ has been a widely studied topic in graph theory,
especially in its form known as the \textit{decycling set problem}
or the \textit{feedback vertex set problem}. A decycling set of a
graph is a subset of its vertices whose deletion yields an acyclic
graph. From this definition, we deduce that a set $S \subseteq V$
is such that $G[S]$ is an induced forest of $G=(V,E)$ if and only
if $V \setminus S$ is a decycling set of $G$. Therefore, finding a
lower bound for $\tau(G)$, the maximum number of vertices in an
induced forest of $G$, amounts to finding an upper bound for
$\phi(G)$, the minimum cardinality of a decycling set of $G$.

Historically, the problem of obtaining an acyclic subgraph of a
graph $G$ by removing vertices was already considered by Kirchhoff
in his work on spanning trees \cite{kirchhoff}. Erd\"{o}s et.\
al.\ also worked on this problem stated in terms of maximum
induced trees in a graph \cite{indtree}. However, finding a
decycling set of a given size in a graph is inherently difficult.
Indeed, this problem has been shown to be NP-complete
\cite{karp1}, even for special families of graphs such as
bipartite graphs, planar graphs or perfect graphs.

On the other hand, there exist polynomial algorithms to solve
instances of this problem in cubic graphs \cite{li_liu},
permutation graphs \cite{liang} and interval graphs
\cite{liang_chang}. Also, tighter bounds or even the exact value
of the decycling number have been determined for graphs such as
grids and cubes in \cite{bau_beineke} and \cite{beineke_vandell}.

For random regular graphs with fixed degree $r$, upper and lower
bounds on the size of a minimum decycling set have been obtained
by Bau et al.~in \cite{bau_wormald}. Their strategy relies on the
analysis of a randomised greedy algorithm that generates a
decycling set of a regular graph as it is exposed in the usual
pairing model of random regular graphs.

We investigate induced forests in $r$-regular graphs with large
girth, where $r \geq 3$ is a fixed integer. By the \textit{girth}
of a graph $G$, we mean the length of a shortest cycle in $G$, if
the graph contains a cycle, or infinity, if it is acyclic. We will
extend the method initiated by Lauer and the second author to find
lower bounds on the size of largest independent sets
\cite{lauer_wormald} in such graphs. The proof involves analysing
the performance of an iterative randomised algorithm that
generates an independent set in a graph. Although their algorithm
is applicable to any graph, the number of iterations allowed is
bounded by a function that increases with the girth, and, because
of this, better bounds can be obtained as the girth increases. We
shall use a similar approach to obtain bounds on the size of an
induced forest in a graph whose girth is large.

More precisely, we will prove the following.
\begin{thm} \label{maintheorem}
Let $\delta>0$ and $r \in \mathbb{N}$. Then, there exists $g>0$
such that every $r$-regular graph $G$ on $n$ vertices with girth
greater than or equal to $g$ satisfies $\tau(G) \geq (\xi(r) -
\delta) n$, where the constants $\xi(r)$ are derived from the
solution of a system of differential equations. Numerical values
are given in the table below for some values of $r$.

\begin{center}
\begin{tabular}{|c|c|c|}
\hline $r$ & $\xi(r)$ & $\Xi(r)$\\
\hline 3 & 0.7268  & 0.2732\\
\hline 4 & 0.6045  & 0.3955\\
\hline 5 & 0.5269  & 0.4731\\
\hline 6 & 0.4711  & 0.5289\\
\hline 7 & 0.4283  & 0.5717\\
\hline 8 & 0.3940  & 0.6060\\
\hline 9 & 0.3658  & 0.6342\\
\hline 10 & 0.3419 & 0.6581\\
\hline
\end{tabular}

\vspace{1em} Table~1. Lower bounds on $\tau(G)$ and, in the last
column, upper bounds on $\phi(G)$, where $G$ is an $r$-regular
graph with sufficiently large girth.
\end{center}

\end{thm}
An actual formula for the constants $\xi(r)$ will be given in
Section \ref{finalsection}.

Consider $\Xi(r)=1-\xi(r)$ in the above table. Then, for any fixed
$r \geq 3$ and $\delta>0$, $(\Xi(r)+\delta)n$ gives an upper bound
on the number of vertices in a minimum decycling set of an
$r$-regular graph $G$ with girth greater than or equal to the
positive integer $g$ referred to in the theorem. For all values of
$r$ tested, with the exception of $r=3$, these are the best bounds
known for regular graphs with large girth. For $r=3$, it can be
shown that $\tau(G)=0.75$ for every $3$-regular graph with
sufficiently large girth as a consequence of a result on
fragmentability of graphs \cite{edwards_farr}.

Observe that, if $G$ is a graph with maximum degree $r$, then we
can create an $r$-regular graph by taking copies of $G$ and
joining some pairs of vertices from different copies so as to make
the resulting graph $G^\prime$ $r$-regular. This can be done
without decreasing the girth if sufficiently many copies of $G$
are used. Moreover, we have the inequality $\tau(G) \geq
\tau(G^\prime)$ because the copy of $G$ containing the most
vertices in a largest induced forest in $G^\prime$ satisfies this
property. Thus, the theorem immediately implies the following
result.

\begin{cor}
Let $\delta>0$ and $r \in \mathbb{N}$. Then, there exists $g>0$
such that every graph $G$ on $n$ vertices with maximum degree $r$
and girth greater than or equal to $g$ satisfies $\tau(G) \geq
(\xi(r) - \delta) n$, where $\xi(r)$ is given in Table~1 for some
values of $r$.
\end{cor}

Furthermore, the second author~\cite{wormald1} (see
also~\cite{Wthesis}) and Bollob\'{a}s~\cite{bollobas1}
independently proved results implying that, if $G$ is a random
$r$-regular graph on $n$ vertices, and $g$ is any positive
integer, then $G$ asymptotically almost surely has $o(n)$ cycles
of length at most $g$. (For a sequence of probability spaces
$\Omega_n$, $n \geq 1$, an event $A_n$ of $\Omega_n$ occurs
asymptotically almost surely, or a.a.s.\ for brevity, if $\lim_{n
\rightarrow \infty} \pr(A_n)=1$.) So, $G$ a.a.s.\ can be turned
into a graph $G^{\prime}$ with maximum degree $r$ and girth at
least $g$ by deleting $o(n)$ of its vertices. By the previous
corollary, given $\delta>0$, we can find $g>0$ such that
$G^{\prime}$ contains an induced forest with at least $(\xi(r) -
\delta) \big( n - o(n) \big)$ vertices. If we delete from the
forest all the vertices adjacent to a vertex of $V(G) \setminus
V(G^{\prime})$, we have an induced forest of $G$ with at least
$(\xi(r) - \delta) \big( n-o(n) \big) - o( n)$ vertices. This
leads to the following result.
\begin{cor}
Let $r \in \mathbb{N}$ and let $G$ be a random $r$-regular graph
on $n$ vertices. Fix $\epsilon >0$. Then a.a.s.\ $G$ contains an
induced forest with $(\xi(r)-\epsilon)n$ vertices, where $\xi(r)$
is the constant given in Table~1.
\end{cor}
For $r \geq 4$, these bounds coincide with the  corrected version of 
the best
bounds known for random regular graphs obtained in
\cite{bau_wormald}. The need for correction arises from a
fault in the latter part of the argument, which relied upon greedily 
growing an
induced forest in a random regular graph. The number of
uninvestigated edges leading out of the forest at time $t$ was
denoted $Y(t)$. Differential equations were set up which describe
the likely behaviour of $Y(t)$, and it was shown that the actual
behaviour is close to the likely behaviour a.a.s.  The argument is
valid as long as $Y(t)>0$. Unfortunately, the equations were
traced after that point, so the bounds quoted
in~\cite{bau_wormald} for $r\ge 4$ are invalid. It is quite easy
to correct this. At the time $Y(t)$ falls to 0, all vertices
adjacent to the growing forest are adjacent to at least two
vertices in the forest and cannot be added to it without creating
a cycle, and the forest is actually a tree $T$. From the values of
the variables at that point, it is easy to obtain the number of
vertices not adjacent to any vertices in the tree, and also the
number of edges in the subgraph $H$ induced by them. Since the
average degree is less than and bounded away from 1, and $H$ is
uniformly distributed, given its degree sequence, it is easy to
show that the number of vertices in cycles of $H$ is a.a.s.\ small
(say less than $\log n$), and thus almost all vertices of $H$ can
be added to $T$ to obtain an induced forest in the original graph.
This argument gives, for the random graph, the same bounds as in
Table~1.

The main goal of this paper is to establish
Theorem~\ref{maintheorem}, whose proof is structured as follows.
We first introduce a randomised greedy algorithm that finds an
induced forest of a given graph. As with the discussion above on
random regular graphs, the final part of this algorithm adds
almost all the vertices in a set of ``leftover" vertices. When
this algorithm is applied to an $r$-regular graph $G$ with
sufficiently large girth, its expected performance leads to the
bounds in Table~1, and hence guarantees the existence of an
induced forest on the same proportion of vertices by the first
moment principle. To estimate the expected performance of the
algorithm, we shall establish preliminary lemmas that help us
understand the behaviour of our algorithm, which will then be used
to derive a system of recurrence equations involving the
cardinality of the set of vertices in the induced forest. Finally,
we shall approximate this system of recurrence equations by a
system of ordinary differential equations whose solution provides
us with the bounds mentioned above.

Our method also produces (weaker) bounds on $\tau(G)$ if a
specific lower bound on the girth of $G$ is given. However, we do
not compute the precise constants for any particular bound on
girth.

We remark here that the method initiated in~\cite{lauer_wormald}
and developed further in this paper can clearly be applied and
adapted to obtain a wide range of results on large sets of
vertices or edges in bounded degree graphs with large girth. A
particularly powerful extension which the authors are planning is
to permit prioritisation of a number of alternative steps in the
greedy algorithm. Such steps are used in the most powerful
algorithms known for finding independent sets or dominating sets
in random regular graphs; see~\cite{W} and~\cite{DW}.

\section{ An algorithm}

We introduce an algorithm that will help us find a large induced
forest in a graph $G=(V,E)$. At any given step of the algorithm,
we shall associate colours with the vertices of the graph as
follows. The colour purple is assigned to vertices in a set $P$, a
subset of $V$ that induces a subgraph of $G$ with ``a few" cycles
only. A vertex is blue if it is not yet in $P$, but could join it
in the next iteration, whereas orange is assigned to vertices
whose addition to $P$ would yield cycles in $G[P]$. The remaining
vertices are coloured white and are the vertices not adjacent to
vertices of the forest.

\begin{alg}\label{alg1}

\

\noindent Input: A graph $G$, a positive integer $N$ and a pair of
probabilities $(p_0,\pix )$.

\begin{enumerate}
\item[1.] Start with all the vertices of the graph coloured
white. In the first step, colour each vertex purple with probability
$p_0$, at random, independently of all others. Non-purple vertices
are coloured blue if they have exactly one purple neighbour and
orange if they have at least two purple neighbours.

\item[2.] At each step $i$, choose blue vertices randomly and
independently with probability $\pix$ and colour them purple.
The sets of blue and orange vertices are updated using the rule given 
in 1. We refer to the set of white vertices as $W$, to
the set of blue vertices as $B$ and to the set of purple vertices
as $P$. Repeat this iteratively for $N$ steps.

\item[3.] Create a set $\bar{P} \subseteq P$ by deleting any pair of
adjacent vertices added to $P$ in a same step.

\end{enumerate}

\noindent Output: The acyclic set $\bar{P}$ and the set of white
vertices $W$.

\end{alg}

In the first phase, the roots of the induced trees are chosen and
coloured purple, and vertices that could be added to the trees
without creating cycles or connecting distinct components are
coloured blue. In each step of the second phase, the forest is
extended by choosing blue vertices and adding them to $P$, and 
at each step the
colours associated with each vertex are updated so that
the sets of white, blue and orange vertices at the end of each
step represent the vertices with 0, 1, and more than one, purple
neighbours, respectively. Note that it would be possible to
alter $p$ at each step, and this would be useful if optimising the
algorithm for the set of graphs with particular girth (as done 
in~\cite{lauer_wormald} for independent sets), but we do
not do this here.

The graph $G[P]$ at the end of Phase 2 is not
necessarily acyclic. As a matter of fact, it may happen that two
neighbouring blue vertices are added to the forest in the same
step and create a cycle. So, the set of purple vertices induces a
subgraph with ``a few" cycles, and the set of orange vertices is
``almost" a decycling set of the graph. These cycles are broken in
the third phase of the algorithm.

A drawback to the analysis of Algorithm~\ref{alg1} in its original
version is that the random selection of vertices at a given step
depends on the outcome of the previous steps. To avoid this, we
introduce an equivalent model for which the random choices are
uniform over the whole set of vertices. Indeed, with each vertex
$v \in V$, we shall first associate a random sequence of labels
$S(v) \subseteq \{0,1,2,\ldots\}$ so that label $i$ is in $S(v)$
independently at random with probability $p_0$, if $i=0$, or
$\pix$, if $i \geq 1$. In other words, we choose sets of vertices
at times 0, 1, $\ldots$, and assign to a vertex $v$ the labels
$\{i:~v \textrm{ was chosen at time } i\}$. In the context of our
algorithm, we shall then consider the set of vertices with label
$0$ to be the set of vertices selected in Phase 1 and use vertices
with label $i$ to recreate the set of vertices added to $P$ at
step $i$ in Phase 2 of our algorithm. It is clear that some of the
labels are ill-suited. For instance, a vertex with label 1 will
not be selected to join $P$ at step 1 if it also has label 0, in
which case it already belongs to $P$, or if none of its neighbours
has label 0, in which case it is not blue after the first phase of
the algorithm. This motivates a classification of the labels as
relevant or irrelevant, that is, as labels that represent an
action of our algorithm or as labels that do not.

\begin{defn}{\textit{Relevant and irrelevant labels}}

\noindent Let $G=(V,E)$ be a graph, and, for every $v \in V$, let
$S(v) \subseteq \mathbb{N}$ be the set of labels associated with
$v$. We define relevant labels inductively (labels that are not
relevant are said to be irrelevant). A label $i$ is relevant for
$v$ if:

\noindent I. $i=0 \in S(v)$, or

\noindent II. $i \in S(v)$, $j$ is irrelevant for $v$ for all $j <
i$, and there is a unique neighbour of $v$ with a relevant label
strictly smaller than $i$.

\end{defn}

The sets of vertices with relevant label equal to $i$ are denoted
by $R_i$, while the ones with relevant label less than or equal to
$i$ are denoted by and $R_{\leq i}$. We refer to the sequence
$[S(v): v \in V]$ as $\mathcal{S}$. Now, for each $l \in \mathbb{N}$, 
the sequence $\mathcal{S}$ may
be used to construct a colouring of $G$ with colours purple, blue,
white and orange.

\begin{defn}{\textit{Colouring of $G$ at time $l$}}

\noindent Given a graph $G$ and a sequence $\mathcal{S}$ as above,
the colouring of $G$ at time $l \in \mathbb{N}$ is the function
assigning colours purple, blue, orange and white to the vertices
of $G$ defined as follows. Given $u \in V$,

\begin{enumerate}

\item[(a)] $u$ is white if $u \notin R_{\leq l}$ and $v \notin
R_{\leq l}$, for all $v \in N(u)$, where $N(u)$ denotes the
neighbourhood of $u$.

\item[(b)] $u$ is blue if $u \notin R_{\leq l}$ and there is a
unique $v \in N(u)$ such that $v \in R_{\leq l}$.

\item[(c)] $u$ is orange if $u \notin R_{\leq l}$ and there
exist distinct $v,w \in N(u)$ with $v,w \in R_{\leq l}$.

\item[(d)] $u$ is purple if $u \in R_{\leq l}$.

\end{enumerate}

\end{defn}

It is clear from this definition that the colouring of $G$ at time
$l$ is fully determined by the sequence $[S(v) \cap \{0,\ldots,l\}
~:~ v \in V]$. Moreover, this colouring coincides with the
colouring of the graph induced by our algorithm if we assume the
set $P$ after $k$ steps to be $R_{\leq k}$, as formalised by the
next lemma.

\begin{lemma}
Let $G=(V,E)$ be a graph, and consider a subgraph $H$ of $G$ and a
colouring $c$ of $H$ with colours purple, blue, orange and white.
Then, the following events have the same probability:

\begin{itemize}
\item[(i)] the colouring of $G$ at time $l$ induced by the sequence
$\mathcal{S}=[S(v) ~:~ v \in V(G)]$ restricted to $H$ is equal to
$c$, where $\mathcal{S}$ is obtained by adding each nonnegative
integer $i$ to $S(v)$ independently with probability $p_0$, if
$i=0$, or $\pix$, if $i \geq 1$, for all $v \in V$.

\item[(ii)] Algorithm~\ref{alg1} applied to $G$ obtains $c$ as the
colouring of $H$ after step $l$.
\end{itemize}
\end{lemma}

\proof We modify Phase 2 our algorithm to allow all vertices to be
chosen uniformly at random with probability $\pix$, instead of
restricting our choices to blue vertices. However, no action is
taken if a non-blue vertex is selected. So, these extra ``dummy"
choices do not alter the probability of a given colouring of $G$
being obtained and our result follows. \qed

In the remainder of this paper, we shall work in the probability
space of the sequence $\mathcal{S}$ of sets of labels. So, each
time a colouring of graph $G$ is mentioned, the colouring induced
by $\mathcal{S}$ is meant.

\section{Independence lemmas}

We prove results that allow us to compute the probability, using
local information only, of a vertex of an $r$-regular graph $G$
being assigned some given colour at time $i$. Henceforth, we shall
fix an $r$-regular graph $G = (V, E)$ with girth $g$ and consider
a sequence of sets $\mathcal{S}=[S(v) ~:~ v \in V(G)]$, where $i
\in \mathbb{N}$ is in $S(v)$ with probability $p_0$, if $i=0$, or
$\pix$, if $i \geq 1$, for all $v$.

\begin{lemma}\label{lem0}
Let $G=(V,E)$ be a graph and consider a sequence of sets of labels
$\mathcal{S} = [S(v) ~:~ v \in V]$. Given $u \in V$, define a
sequence of sets of labels $\mathcal{S}^{\prime}$ by replacing, in
$\mathcal{S}$, $S(u)$ by some set $S^{\prime}(u)$. Let $w$ be a
vertex of $G$ whose colours at time $i$ with respect to
$\mathcal{S}$ and $\mathcal{S}^{\prime}$ differ, where $i$ is a
nonnegative integer.

Then, there exists a path $\mathcal{P}$ from $u$ to $w$ for which
every vertex except possibly $w$ gained or lost a relevant label
less than or equal to $i$ when $\mathcal{S}$ was replaced by
$\mathcal{S}^{\prime}$. Moreover, the relevant labels gained or
lost by each vertex along the path are in strictly increasing
order when the path is considered from $u$ to $w$.
\end{lemma}

\proof The proof is by induction on $i$. For $i=0$, since the
colour of $w$ at time $0$ has changed after replacing $S(u)$ by
$S^{\prime}(u)$, it must be that $u$ has gained or lost relevant
label $0$ and that either $u=w$ or $u$ and $w$ are neighbours. In
both cases, $\mathcal{P}=(u,w)$ satisfies the conditions in the
statement of this lemma.

Now, let $i>0$ and assume that this result holds earlier. If $u=w$
nothing needs to be done, so suppose that this is not the case.
Since the colour of $w$ changed at time $i$, there exists a
neighbour $w^{\prime}$ of $w$ that gained or lost a relevant label
smaller than or equal to $i$. If $w^{\prime}=u$, our result is
clearly true, so suppose that they are distinct. Then, the
relevant label gained or lost by $w^{\prime}$ is not equal to 0
and, by the definition of relevant label, there is a neighbour
$w^{\prime \prime}$ of $w^{\prime}$ that gained or lost a relevant
label at a time $j$ strictly smaller than the relevant label
gained or lost by $w^{\prime}$. In particular, the colour of
$w^{\prime \prime}$ changed at time $j$, so, by induction, there
is a path $\mathcal{P}^{\prime \prime}$ from $u$ to $w^{\prime
\prime}$ under the conditions of the lemma. Thus, the path
$\mathcal{P}$ obtained by appending vertices $w^{\prime}$ and $u$
to $\mathcal{P}^{\prime \prime}$ satisfies the required
properties. \qed

\begin{cor}\label{lem1}
Let $u \in V(G)$ and $i \in \mathbb{N}$. Then, for any given
colour $c$ and any collection of subsets $S^{\prime}_v$ of
$\mathbb{N}$, where $v$ ranges over the vertices at distance at
least $i+2$ of $u$, the event that $u$ has colour $c$ at time $i$
is independent of the event that $S(v)=S^{\prime}_v$.
\end{cor}

\proof  It is sufficient to show that, if
$\hat{\mathcal{S}}=[\hat{S}(v) ~:~ v \in V(G)]$ is any given
family of sets of labels and new sets $S^{\prime}(v)$ are assigned
to each vertex $v$ satisfying $d(u,v) \geq i+2$, then the colour
of $u$ at time $i$ relative to $\hat{\mathcal{S}}$ is the same as
the colour of $u$ at time $i$ relative to $\mathcal{S}^{\prime}$,
where $\mathcal{S}^{\prime}$ is obtained by replacing each
$\hat{S}(v)$ by $S^{\prime}(v)$.

We now prove this sufficient condition. Suppose for a
contradiction that the colours of $u$ with respect to
$\hat{\mathcal{S}}$ and $\mathcal{S}^{\prime}$ differ, and order
the vertices $v \in V$ satisfying $d(u,v) \geq i+2$ as $v_1, v_2,
\ldots ,v_m$. Consider, for $l \in \{0,\ldots,m\}$, the sequences
$\mathcal{S}_l$ obtained from $\hat{\mathcal{S}}$ by replacing
$\hat{S}(v_1),\ldots,\hat{S}(v_l)$ by
$S^{\prime}(v_1),\ldots,S^{\prime}(v_l)$. Our assumption implies
the existence of $j$ such that the colours of $u$ with respect to
$\mathcal{S}_j$ and $\mathcal{S}_{j+1}$ are distinct. By
Lemma~\ref{lem0}, there is a path $\mathcal{P}$ in $G$ from
$v_{j+1}$ to $u$ such that every vertex except possibly $u$ gained
or lost a relevant label less than or equal to $i$ when
$\mathcal{S}_j$ was replaced by $\mathcal{S}_{j+1}$. Also, the
relevant labels gained or lost on each vertex along the path are
in strictly increasing order when the path is considered from
$v_{j+1}$ to $u$. In particular, $\mathcal{P}$ contains at most
$i+2$ vertices, i.e., $d(u,v_{j+1}) \leq i+1$, a contradiction.
\qed

Let $B_i$ and $W_i$ denote the sets of vertices coloured blue and
white at time $i$, respectively.

\begin{cor}\label{lem3}
Let $u \in V$ and let v be one its neighbours. Then, the
probabilities $\pr(u \in W_i)$, $\pr(u \in B_i)$,
$\pr(u \in W_i \wedge v \in W_i)$, $\pr(u \in B_i
\wedge u \in W_i)$ and $\pr(u \in B_i \wedge v \in B_i)$
are independent of $u$ and $v$ whenever $2i < g-3$. Moreover, if
we let $w$ be a neighbour of $u$ distinct from $v$, $\pr(v
\in B_i \wedge u \in B_i \wedge w \in R_{\leq i})$ does not depend
on $u$, $v$ or $w$.
\end{cor}

\proof We know from Corollary~\ref{lem1} that the colour of $u$ at
time $i$ depends only on the sets of labels of vertices at
distance at most $i+1$ from $u$. In other words, $u$ is fully
determined by the sets of labels in the subgraph $G_u=G[\{v :
d_{G}(u,v) \leq i+1\}]$. But our restriction on $i$ implies that,
for every $u \in V$, the graphs $G_u$ are isomorphic. Our first
two claims immediately follow, since distinct vertices are
assigned sets of labels independently with the same probability.
It is clear that an analogous argument can be used to prove the
remaining statements. \qed

\begin{lemma}\label{lem2}
Let $u \in V$ and let $v_1,\ldots,v_r$ be its neighbours. Fix $i,k
\in \mathbb{N}$ such that $2(i+2)<g-2$ and consider, for each $j
\in \{1,\ldots,r\}$, the tree $T_j^{2}$ rooted at $v_j$ given by
the component of $G[\{v : d_{G}(u,v) \leq 2\} - u]$ containing
$v_j$.

Then, the following assertions hold.

\noindent 1. Let  $X_1,\ldots,X_r$  be colourings of the tree
isomorphic to the rooted trees $T_j^{2}$ (the isomorphism is a
consequence of our restriction on $i,k$). Then, conditional upon
$u \in W_i$, the events $E_1,\ldots,E_r$ are mutually independent,
where $E_j$ stands for the event that $T_j^{2}$ has colouring
$X_j$ at time $i$.

\noindent 2. Conditional upon $u \in B_i$ and $v_l \in R_{\leq i}$
for some $\ \in \{1,\ldots,r\}$, the same events $E_j$ are
mutually independent for all $j \neq l$.
\end{lemma}

\proof By Lemma~\ref{lem0}, the colour of a vertex $w$ at time $i$
is altered when replacing $S_1(v) \times \cdots \times S_r(v)$ by
$S^{\prime}_1(v) \times \cdots \times S^{\prime}_r(v)$ only if
there is a path $\mathcal{P}$ from $v$ to $w$ such that all the
vertices on $\mathcal{P}$ that are not purple at time $i$ with
respect to $\mathcal{S}$ have a different colour with respect to
$\mathcal{S}^{\prime}$. This is because, given any non-purple
vertex $w^{\prime}$ at time $i$ lying on $\mathcal{P}$, it either
gains a relevant label, in which case it is purple at time $i$
with respect to $\mathcal{S}^{\prime}$, or it is equal to $w$, in
which case its colour changes by assumption.

Now, if $w$ and $v$ are vertices in different branches with
respect to $u$, our restriction on $i$ implies by
Corollary~\ref{lem1} that any path from $v$ to $w$ that is short
enough for every interior vertex to gain or lose a relevant label
passes through $u$. Hence, conditional upon $u$ being white,
changes in $S_1(v)\times \cdots \times S_r(v)$ do not affect the
colour of $w$ at time $i$.

Moreover, we also know by Corollary~\ref{lem1} that the colour at
time $i$ of vertices at distance at most two from $u$ are not
affected by changes in the set of labels of vertices whose
distance to $u$ is greater that $i+3$. Let $V_{u,i,k}$ be the set
of vertices in $G$ at distance at most $k+i+1$ from $u$, excluding
vertex $u$.

So,
\begin{equation*}
\pr(E_1 \wedge E_2 \wedge \ldots \wedge E_r ~|~ u \in W_i) =
\sum_{\star} \pr(S(v)=S_v, \forall v \in V_{u,i,k}~|~ u \in W_i),
\end{equation*}
where $\sum_{\star}$ denotes the sum over vectors $(S_v~:~v \in
V_{u,i,k})$ such that the event $S(v)=S_v,~\forall v \in
V_{u,i,k}$, implies $E_1 \wedge E_2 \wedge \ldots \wedge E_r$.
Now, observe that our restriction on $i,k$ implies that the trees
$T_j^{(k+i+1)}$ are all disjoint. In particular, we can first sum
over sets of labels of vertices in $T_1^{(k+i+1)}$ (notation
$\sum_{\star \star}$) and then over the remaining vertices
(notation $\sum_{\star \star \star}$) to obtain
\begin{equation*}
\pr(E_1 \wedge E_2 \wedge \ldots \wedge E_r ~|~ u \in W_i) =
\sum_{\star \star} \sum_{\star \star \star} \pr(S(v)=S_v, \forall
v \in V_{u,i,k} ~|~ u \in W_i).
\end{equation*}
Using conditional probability and rearranging the sum, this
becomes
\begin{equation*}
\begin{split}
& \sum_{\star \star} \pr(S(v)=S_v, \forall v \in
 T_1^{(k+i+1)} ~|~ u \in W_i) \times \\
& ~~~ \times \sum_{\star \star \star} \pr(S(v)=S_v, \forall
v \notin T_1^{(k+i+1)} ~|~ (u \in W_i) \wedge (S(v)=S_v, \forall v
\in
T_1^{(k+i+1)}))\\
& = \sum_{\star \star} \pr(S(v)=S_v, \forall v \in T_1^{(k+i+1)}
~|~ u \in W_i) ~\pr(E_2 \wedge \ldots \wedge E_r
~|~ u \in W_i)\\
&= \pr(E_1 ~|~ u \in W_i)~\pr(E_2 \wedge \ldots \wedge E_r ~|~ u
\in W_i).
\end{split}
\end{equation*}
These manipulations can be done since, conditional upon $u$ being
white, changes in $S(v)$ do not affect the colours of other
branches, for any $v \in T_1^{(k+i+1)}$.

Repeating this argument for the remaining branches, we obtain
\begin{equation*}
\pr(E_1 \wedge E_2 \wedge \ldots \wedge E_r ~|~ u \in W_i) =
\prod_{j = 1}^r \pr(E_j ~|~ u \in W_i),
\end{equation*}
and our first claim is true.

For the second part, we proceed analogously by leaving both the
blue vertex $u$ and the branch of its neighbour with relevant
label untouched, and then summing over all possibilities of labels
for vertices in the other branches. \qed

\section{Applications of the Independence Lemmas}

In this section, the independence results of the previous section will 
be used to obtain
recurrence equations relating the probabilities of events that are
important in the analysis of Algorithm \ref{alg1}. We introduce some
notation. Let $u$ be a vertex of graph $G$. An arbitrary neighbour
of $u$ will be denoted by $v$, while we use $v_1,\ldots,v_r$ to
refer to the set of neighbours of $u$. When $u$ has a neighbour
with relevant label, this will be referred as $v_k$ and we shall
assume that $v \neq v_k$.

Furthermore, for any $i \geq 0$, we know by Corollary~\ref{lem3}
that the quantities $w_i = \pr(u \in W_i)$, $b_i =
\pr(u \in B_i)$, $q_i = \pr(u \in W_i \wedge v \in
W_i)$, $s_i = \pr(u \in B_i \wedge v \in W_i)$ and $t_i =
\pr(u \in B_i \wedge v \in B_i)$, or even $\pr(u \in
B_i \wedge v \in W_i \wedge v_k \in R_{\leq i})$ and $\pr(u
\in B_i \wedge v \in B_i \wedge v_k \in R_{\leq i})$, do not
depend on $u$, $v$ or $k$. We now let $i \geq 1$ and establish the
following consequences of the previous independence lemmas.

\begin{cor}\label{cor1}

\

\begin{itemize}
\item[(i)] Let $J=\{j_1,\ldots,j_k\} \subseteq \{1,\ldots,r\}$. Then,
\begin{equation*}
\pr(v_j \notin R_i, ~\forall j \in J ~|~ u \in
W_{i-1})=\prod_{j \in J} \pr(v_j \notin R_i ~|~ u \in
W_{i-1}).
\end{equation*}
\item[(ii)] Let $J \subseteq \{1,\ldots,r\} \setminus \{k\}$. Then,
\begin{equation*}
\begin{split}
\pr(v_j \notin R_i, ~\forall j \in J ~|~ u \in B_{i-1} \wedge v_k&
\in R_{\leq i-1})= \\
& \prod_{j \in J} \pr(v_j \notin R_i ~|~ u \in B_{i-1}\wedge v_k
\in R_{\leq i-1}).
\end{split}
\end{equation*}
\end{itemize}
\end{cor}

\proof We prove part (i) by induction on $k$. For $k=1$, the
result follows immediately, so let $k>1$ and assume the result
holds for any smaller set $J$.

First observe that, because a vertex receives relevant label $i
\geq 1$ only if it is blue at time $i-1$, it is important to
consider the set of blue neighbours of $u$ at time $i-1$. In light
of this, we associate a vector $\omega \in \mathbb{Z}_2^{k}$ with
the set of neighbours $v_{j_t}$ of $u$ so that $\omega(t)=1$ if
and only if $v_{j_t} \in B_{i-1}$.

Note that, for a vertex not to become purple at time $i$, it
either was not blue at the previous step or it was blue, but $i$
is not contained in its set of labels. Thus,
\begin{align*}
\pr&(v_j \notin R_i, ~\forall j \in J ~|~ u \in W_{i-1})\\
& = \sum_{\omega \in \mathbb{Z}_2^{k}} \pr((v_{j_t} \in B_{i-1}
\wedge i \notin S(v_{j_t}), \forall t \textrm{ with }
\omega(t)=1) ~\wedge\\
& \ \ \ \ \ \ \ \ \ \ \ \ \ \ \  \wedge (v_{j_t} \notin B_{i-1},
\forall t \textrm{ with } \omega(t)=0)~|~u \in W_{i-1}).
\end{align*}
The fact that $S(v)$ contains any nonnegative integer independently at
random (and label $i$ does not influence the colouring at time
$i-1$), together with Lemma~\ref{lem2}, ensures that the events of
the form $(v_{j_t} \in B_{i-1}\wedge i \notin S(v_{j_t}))$ and
$v_{j_u} \notin B_{i-1}$ are mutually independent conditional upon
$u$ being white. So, the equation becomes
\begin{align*}
\pr&(v_j \notin R_i, ~\forall j \in J ~|~ u \in W_{i-1})\\
& = \sum_{\omega \in \mathbb{Z}_2^{k}} \prod_{\{j_t :
\omega(t)=1\}} \pr(v_{j_t} \in B_{i-1} \wedge i \notin S(v_{j_t})
~|~u \in W_{i-1}) ~\times \\
& \ \ \ \ \ \ \ \ \ \ \ \ \ \ \ \times \prod_{\{j_t :
\omega(t)=0\}}
\pr(v_{j_t} \notin B_{i-1}~|~u \in W_{i-1}).\\
& = \pr(v_{j_k} \in B_{i-1} \wedge i \notin S(v_{j_k}) ~|~u \in
W_{i-1}) \sum_{\omega^{\prime} \in \mathbb{Z}_2^{k-1}}
\prod_{\{j_t :
\omega^{\prime}(t)=0\}} \pr(v_{j_t} \notin B_{i-1}~|~u \in 
W_{i-1})~\times \\
& \ \ \ \ \ \ \ \ \ \ \ \ \ \ \ \times \prod_{\{j_t :
\omega^{\prime}(t)=1\}} \pr(v_{j_t} \in B_{i-1} \wedge i \notin
S(v_{j_t}) ~|~u \in W_{i-1})~+\\
& \ \ \ \ \ +\pr(v_{j_k} \notin B_{i-1}~|~u \in
W_{i-1})\sum_{\omega^{\prime} \in \mathbb{Z}_2^{k-1}} \prod_{\{j_t
:
\omega^{\prime}(t)=0\}} \pr(v_{j_t} \notin B_{i-1}~|~u \in 
W_{i-1})~\times \\
& \ \ \ \ \ \ \ \ \ \ \ \ \ \ \ \times \prod_{\{j_t :
\omega^{\prime}(t)=1\}} \pr(v_{j_t} \in B_{i-1} \wedge i \notin
S(v_{j_t}) ~|~u \in W_{i-1})\\
& = \pr(v_{j_k} \in B_{i-1} \wedge i \notin S(v_{j_k}) ~|~u \in
W_{i-1}) \pr(v_j \notin R_i, ~\forall j \in J\setminus \{k\} ~|~ u
\in W_{i-1}) ~+\\
&  \ \ \ \ \ + \pr(v_{j_k} \notin B_{i-1}~|~u \in W_{i-1})\pr(v_j
\notin R_i, ~\forall j \in J\setminus \{k\} ~|~ u \in W_{i-1})
\end{align*}
By induction, this is equal to
\begin{align*}
\pr&(v_{j_k} \in B_{i-1} \wedge i \notin S(v_{j_k}) ~|~u \in
W_{i-1}) \prod_{j \in J\setminus \{k\}} \pr(v_j \notin
R_i ~|~ u \in B_{i-1}\wedge v_k \in R_{\leq i-1}) ~ +\\
&  \ \ \ \ \ + \pr(v_{j_k} \notin B_{i-1}~|~u \in W_{i-1})
\prod_{j \in J\setminus \{k\}} \pr(v_j \notin
R_i ~|~ u \in B_{i-1}\wedge v_k \in R_{\leq i-1})\\
& = \prod_{j \in J} \pr(v_j \notin R_i ~|~ u \in W_{i-1}),
\end{align*}
as required for (i).

An analogous argument gives (ii). \qed

\begin{rem}
This corollary can also be extended to conditioning upon $u \in
W_{i-1} \wedge v \in W_{i-1}$, where $u,v$ are neighbours in $G$
(or any other combination of restrictions on $u,v$ being white or
blue). As a matter of fact, if $u_1,\ldots,u_{r-1}$,
$v_1,\ldots,v_{r-1}$ denote the neighbours of $u,v$ distinct from
$u$ and $v$, and $J,K \subseteq \{1,\ldots,r-1\}$, then
\begin{equation*}
\begin{split}
\pr&((u_j \notin R_i, ~\forall j \in J) \wedge (v_k \notin
R_i, ~\forall k \in K) ~|~ u \in W_{i-1} \wedge v \in W_{i-1})\\
& = \prod_{j \in J} \pr(u_j \notin R_i ~|~ u \in W_{i-1})\prod_{k
\in K} \pr(v_k \notin R_i ~|~ v \in W_{i-1}).
\end{split}
\end{equation*}
This can be obtained by expanding the initial probability into a
sum over vectors $\omega \in \mathbb{Z}_2^{|J|+|K|}$ and then
using the fact that, for any event $E$, we have
$$\pr(E~|~u \in W_{i-1} \wedge v \in W_{i-1}) =
\frac{\pr(E \wedge u \in W_{i-1} ~|~  v \in W_{i-1})}{\pr(u \in
W_{i-1} ~|~ v \in W_{i-1})},$$ so that Lemma~\ref{lem2} can be
applied first with respect to $u \in W_{i-1}$ and then with
respect to $v \in W_{i-1}$. It is clear that similar results can
be stated by conditioning upon other combinations of $u$ and $v$
being white or blue.
\end{rem}

\begin{cor}\label{cor2}
\ 
\begin{itemize}

\item[(i)] $\displaystyle{\pr(u \in W_i ~|~ u \in W_{i-1}) =
\left(1 - \frac{\pix s_{i-1}}{w_{i-1}}\right)^r}$

\item[(ii)] $\displaystyle{\pr(u \in B_i ~|~ u \in W_{i-1}) =
\frac{r \pix s_{i-1}}{w_{i-1}}\left(1 - \frac{\pix
s_{i-1}}{w_{i-1}}\right)^{r-1}}$

\end{itemize}

\end{cor}

\proof For (i), we just observe that, 
for $u$ to cease to be white at time
$i$, at least one of its neighbours has relevant neighbour $i$. Thus,
$$\pr(u \in W_i ~|~ u \in W_{i-1}) = \pr(v_j \notin R_i, ~\forall
j ~|~ u \in W_{i-1}).$$ Now, by Corollary~\ref{cor1}, part (i),
this last expression is equal to
$$\prod_{j=1}^r \pr(v_j \notin R_i ~|~ u \in W_{i-1}).$$
Finally, Corollary~\ref{lem3} guarantees that the probability of
$v_j$ having relevant label $i$ is independent of $v_j$ and equals
the probability of the event that $i \in S(v_j)$ and $v_j$ is
coloured blue at time $i-1$. So,
$$\pr(v_j \in R_i ~|~ u \in W_{i-1}) = \frac{\pix
s_{i-1}}{w_{i-1}},$$ and
$$\pr(u \in W_i ~|~ u \in W_{i-1}) = \left(1 - \frac{\pix
s_{i-1}}{w_{i-1}} \right)^r$$ as a consequence.

Assertion (ii) may be proven using a similar approach. \qed

\begin{cor} \label{cor3}
$\displaystyle{\pr(u \in B_i ~|~ u \in B_{i-1}) = (1-\pix )\left(1
- \frac{r \pix t_{i-1}}{(r-1)b_{i-1}}\right)^{r-1}}$
\end{cor}

\proof The fact that $u$ is blue at step $i-1$ implies that
exactly one of its neighbours $v_1,\ldots,v_r$ has a relevant
label less than or equal to $i-1$. Thus, $$\pr(u \in B_i ~|~ u \in
B_{i-1})=\sum_{k=1}^r \pr(u \in B_i ~|~ u \in B_{i-1} \wedge v_k
\in R_{\leq i-1}) \pr(v_k \in R_{\leq i-1} ~|~ u \in B_{i-1}).$$

Moreover, $u$ remains blue at time $i$ if neither itself nor any
of its neighbours gains a relevant label at time $i$, i.e.,
$$\pr(u \in B_i ~|~ u \in B_{i-1} \wedge v_k \in R_{\leq
i-1}) = (1-\pix )\pr(v_j \notin R_i, ~\forall j \neq k ~|~ u \in
B_{i-1} \wedge v_k \in R_{\leq i-1}).$$ By Corollary~\ref{cor1},
part (ii), we obtain
\begin{align*}
\pr & (v_j \notin R_i, ~\forall j \neq k ~|~ u \in B_{i-1}
\wedge v_k \in R_{\leq i-1})\\
&= \prod_{j \neq k} \pr(v_j \notin R_i ~|~ u \in B_{i-1}
\wedge v_k \in R_{\leq i-1})\\
&= \left(1 - \pix \pr(v_j \in B_{i-1} ~|~ u \in B_{i-1} \wedge v_k
\in R_{\leq i-1})\right)^{r-1}.
\end{align*}
The last equality follows from the fact that $v \in R_i$ only if
it has label $i$ and was blue at time $i-1$.

Finally, we note that
\begin{equation*}
\begin{split}
\pr&(v_j \in B_{i-1} ~|~ u \in B_{i-1} \wedge v_k \in
R_{\leq i-1})= \frac{\pr(v_j \in B_{i-1} \wedge u \in
B_{i-1} \wedge v_k \in R_{\leq i-1})}{\pr(u \in B_{i-1}
\wedge v_k \in R_{\leq i-1})}\\
&= \frac{\pr(v_j \in B_{i-1} \wedge u \in
B_{i-1})\pr(v_k \in R_{\leq i-1} ~|~ v_j \in B_{i-1} \wedge
u \in B_{i-1})}{\pr(u \in B_{i-1})\pr(v_k \in
R_{\leq i-1} ~|~ u \in B_{i-1})}.
\end{split}
\end{equation*}
By Corollary~\ref{lem3}, we conclude that all the neighbours of
$u$ have the same probability of having a relevant label earlier
than the other neighbours, since the probability of having
relevant label $i$ is equal to $p_0$, if $i=0$, or $\pix b_{i-1}$,
if $i \geq 1$, for any vertex. In particular, we must have
$\pr(v_k \in R_{\leq i-1} ~|~ u \in B_{i-1})=\frac{1}{r}$ and
$\pr(v_k \in R_{\leq i-1} ~|~ v_j \in B_{i-1} \wedge u \in
B_{i-1}) = \frac{1}{r-1}$. So,
\begin{equation*}
\begin{split}
\pr&(u \in B_i ~|~ u \in B_{i-1})= (1-\pix ) \sum_{k=1}^r \pr(v_k
\in R_{\leq i-1} ~|~ u \in B_{i-1}) \left(1 -
\frac{r\pix t_{i-1}}{(r-1)b_{i-1}}\right)^{r-1}\\
& = (1-\pix )\left(1 - \frac{r\pix
t_{i-1}}{(r-1)b_{i-1}}\right)^{r-1},
\end{split}
\end{equation*}
with the last equation following from $\sum_{k=1}^r \pr(v_k \in
R_{\leq i-1} ~|~ u \in B_{i-1})=1$. This concludes the proof. \qed
\nick{Punctuation added, format changed:}
\carlosa{Some of the items were followed by commas, some by semi-colons. I replaced the two semi-colons by commas.}
\begin{cor}\label{cor4}

\
\begin{itemize}

\item[(i)] $\displaystyle \pr(u \in W_i \wedge v \in W_i ~|~ u
\in W_{i-1} \wedge v \in W_{i-1}) = \left(1 - \frac{\pix
s_{i-1}}{w_{i-1}}\right)^{2r-2},$
\item[(ii)] $\displaystyle \pr(u \in W_i \wedge v \in B_i ~|~ u
\in W_{i-1} \wedge v \in W_{i-1}) = \frac{(r-1) \pix
s_{i-1}}{w_{i-1}} \left(1 - \frac{\pix
s_{i-1}}{w_{i-1}}\right)^{2r-3},$
\item[(iii)] $\displaystyle \pr(u \in B_i \wedge v \in B_i ~|~ u
  \in W_{i-1} \wedge v \in W_{i-1}) = \frac{(r-1)^2 \pix^2
  s_{i-1}^2}{w_{i-1}^2}\left(1 - \frac{\pix
  s_{i-1}}{w_{i-1}}\right)^{2r-4},$
  \item[(iv)]
$\pr(u \in W_i \wedge v \in B_i ~|~ u \in W_{i-1} \wedge v \in
B_{i-1})$  

$\displaystyle\hspace{3cm}= (1 - \pix ) \left(1 - \frac{\pix
s_{i-1}}{w_{i-1}}\right)^{r-1}  
   \left(1 - \frac{r \pix
t_{i-1}}{(r-1)b_{i-1}}\right)^{r-2},
$
\item[(v)]$\pr(u \in B_i \wedge v \in B_i ~|~ u \in W_{i-1} \wedge v 
\in  
B_{i-1}) $

$\displaystyle\hspace{3cm}= \frac{(r-1) \pix (1-\pix ) 
s_{i-1}}{w_{i-1}}
  \left(1 - \frac{\pix
s_{i-1}}{w_{i-1}}\right)^{r-2}\left(1 - \frac{r \pix
t_{i-1}}{(r-1)b_{i-1}}\right)^{r-2},$
\item[(vi)] $\displaystyle\pr(u \in B_i \wedge v \in B_i ~|~ u
\in B_{i-1} \wedge v \in B_{i-1}) = (1 - \pix )^2 \left(1 -
\frac{r \pix t_{i-1}}{(r-1)b_{i-1}}\right)^{2r-4}.$
\end{itemize}
\end{cor}

\proof Let $u_1,\ldots,u_{r-1}$ be the neighbours of $u$ other
than $v$ and $v_1,\ldots,v_{r-1}$ be the neighbours of $v$
distinct from $u$. Then,
\begin{equation*}
\begin{split}
\pr&(u \in W_i \wedge v \in W_i ~|~ u \in W_{i-1} \wedge v
\in
W_{i-1})\\
& = \pr(u_1,\ldots,u_{r-1},v_1,\ldots,v_{r-1} \notin R_i ~|~ u
\in W_{i-1} \wedge v \in W_{i-1})\\
& = \prod_{k=1}^{r-1} \pr(u_j \notin R_i ~|~ u \in W_{i-1})
\prod_{k=1}^{r-1} \pr(v_j \notin R_i ~|~ v \in W_{i-1})\\
&= \left(1 - \frac{\pix s_{i-1}}{w_{i-1}}\right)^{2r-2}.
\end{split}
\end{equation*}
This is based on the remark after Corollary~\ref{cor1}.

A similar strategy leads to the other formulae. \qed

\section{Differential Equations}

Using the expressions calculated in the last section, we can now
determine recursive formulae for the variables introduced for the
analysis of our algorithm.

\begin{enumerate}

\item[1.] Formula for $w_i$:
\begin{equation*}
\begin{split}
w_i &= \pr(u \in W_i) = \pr(u \in W_i \wedge u \in
W_{i-1})\\
&= \pr(u \in W_{i-1})\pr(u \in W_i ~|~ u \in W_{i-1}) = w_{i-1}
\left(1-\frac{\pix s_{i-1}}{w_{i-1}}\right)^r.
\end{split}
\end{equation*}

\item[2.]  Formula for $b_i$:
\begin{equation*}
\begin{split}
b_i &= \pr(u \in B_i) = \pr(u \in B_i \wedge u \in
B_{i-1}) + \pr(u \in B_i \wedge u \in
W_{i-1})\\
&= \pr(u \in B_{i-1})\pr(u \in B_i ~|~ u \in
B_{i-1}) +  \pr(u \in W_{i-1})\pr(u \in B_i ~|~ u
\in
W_{i-1}) \\
&= b_{i-1}(1-\pix )\left(1-\frac{r \pix
t_{i-1}}{(r-1)b_{i-1}}\right)^{r-1} + r \pix s_{i-1}
\left(1-\frac{\pix s _{i-1}}{w_{i-1}}\right)^{r-1}.
\end{split}
\end{equation*}

\item[3.] Formula for $q_i$:
\begin{equation*}
\begin{split}
q_i &= \pr(u \in W_i \wedge v \in W_i)\\
&= \pr(u \in W_{i-1} \wedge v \in W_{i-1})\pr(u \in
W_i \wedge v \in W_i ~|~ u \in W_{i-1} \wedge
v \in W_{i-1})\\
&= q_{i-1}\left(1 - \frac{\pix s_{i-1}}{w_{i-1}}\right)^{2r-2}.
\end{split}
\end{equation*}

\item[4.] Formula for $s_i$:
\begin{equation*}
\begin{split}
s_i &= \pr(u \in B_i \wedge v \in W_i) \\
&= \pr(u \in B_{i-1} \wedge v \in W_{i-1})\pr(u \in
B_i \wedge v \in W_i ~|~ u \in B_{i-1} \wedge
v \in W_{i-1}) + \\
&+ \pr(u \in W_{i-1} \wedge v \in W_{i-1})\pr(u \in
B_i \wedge v \in W_i ~|~ u \in W_{i-1} \wedge
v \in W_{i-1})\\
&= s_{i-1}(1-\pix )\left(1-\frac{\pix
s_{i-1}}{w_{i-1}}\right)^{r-1}\left(1-\frac{r \pix
t_{i-1}}{(r-1)b_{i-1}}\right)^{r-2}~+\\
& + \frac{(r-1) \pix q_{i-1} s_{i-1}}{w_{i-1}}\left(1-\frac{\pix
s_{i-1}}{w_{i-1}}\right)^{2r-3}.
\end{split}
\end{equation*}

\item[5.] Formula for $t_i$:
\begin{equation*}
\begin{split}
t_i &= \pr(u \in B_i \wedge v \in B_i)\\
&= \pr(u \in B_{i-1} \wedge v \in B_{i-1})\pr(u \in
B_i \wedge v \in B_i ~|~ u \in B_{i-1} \wedge
v \in B_{i-1})~+ \\
&+ \pr(u \in B_{i-1} \wedge v \in W_{i-1})\pr(u \in
B_i \wedge v \in B_i ~|~ u \in B_{i-1} \wedge
v \in W_{i-1})~+ \\
&+ \pr(v \in W_{i-1} \wedge u \in B_{i-1}) \pr(u \in
B_i \wedge v \in B_i ~|~ u \in W_{i-1} \wedge
v \in B_{i-1})~+ \\
&+ \pr(v \in W_{i-1} \wedge u \in W_{i-1}) \pr(u \in
B_i \wedge v \in B_i ~|~ u \in W_{i-1} \wedge
v \in W_{i-1})\\
&= t_{i-1}(1-\pix )^2\left(1-\frac{r \pix
t_{i-1}}{(r-1)b_{i-1}}\right)^{2r-4}
+\\
& + 2 s_{i-1}(1-\pix ) \left(1-\frac{r \pix
t_{i-1}}{(r-1)b_{i-1}}\right)^{r-2}\frac{(r-1)\pix 
s_{i-1}}{w_{i-1}}\left(1-\frac{\pix s_{i-1}}{w_{i-1}}\right)^{r-2} + \\
& + q_{i-1}
\frac{(r-1)^2\pix^2s_{i-1}^2}{w_{i-1}^2}\left(1-\frac{\pix
s_{i-1}}{w_{i-1}}\right)^{2r-4}.
\end{split}
\end{equation*}
\end{enumerate}

We need to evaluate $w_0$, $b_0$, $q_0$, $s_0$ and $t_0$ to have
the necessary initial conditions for solving the system of
recurrence equations found above. It is easy to see that $w_0 =
\pr(u \in W_0) = (1-p_0)^{r+1}$ and $b_0 = rp_0(1-p_0)^r$,
since for the former neither $u$ nor its neighbours can have
relevant label 0, and for the latter $u$ cannot have relevant
label 0, but exactly one of its neighbours must have it.

Now,
\begin{equation*}
q_0 = \pr(u \in W_0 \wedge v \in W_0) = (1-p_0)^{2r},
\end{equation*}
since the event $v \in W_0 \wedge u \in W_0$ is equivalent to
neither $u,v$ nor any of their other neighbours being chosen in
the first phase of the algorithm (and each vertex is chosen
independently with probability $p_0$).

The equation for $s_0$ is given by
\begin{equation*}
s_0 = \pr(u \in B_0 \wedge u \in W_0) =
(r-1)p_0(1-p_0)^{2r-1}
\end{equation*}
because $v \in B_0 \wedge u \in W_0$ occurs when $u,v$ are not
chosen, no neighbours of $u$ are chosen and precisely one
neighbour of $v$ is chosen.

Finally, the equation for $t_0$ is
\begin{equation*} t_0 =
\pr(u \in B_0 \wedge v \in B_0) =
(r-1)^2p_0^2(1-p_0)^{2r-2}
\end{equation*}
with similar justification.

The recurrence equation for $w_i$ obtained at the beginning of
this section can be seen as
$$w_i  =  w_{i-1} - \pix r s_{i-1} + O(\pix^2).$$ For $\pix$ small,
the term $O(\pix^2)$ should only have a minor influence.
Similarly, each of the other equations of the system of recurrence
equations can be rewritten as a main term added to a term of the
order of $\pix^2$. By ignoring the latter, we obtain the following
auxiliary system of recurrence equations:
\begin{equation}\label{modifiedsystem}
\begin{split}
& w^{\prime}_i = w^{\prime}_{i-1} - \pix r s^{\prime}_{i-1}\\
& b^{\prime}_i =  b^{\prime}_{i-1} + \pix \left( - b^{\prime}_{i-1}- r 
t^{\prime}_{i-1} + r s^{\prime}_{i-1} \right)\\
& q^{\prime}_{i} = q^{\prime}_{i-1} -\pix 
\frac{(2r-2)q^{\prime}_{i-1}s^{\prime}_{i-1}}{w^{\prime}_{i-1}}\\
& s^{\prime}_{i} = s^{\prime}_{i-1} + \pix \left( -s^{\prime}_{i-1} + 
\frac{(r-1) q^{\prime}_{i-1} s^{\prime}_{i-1}}{w^{\prime}_{i-1}} 
\right.\\
& \ \ \ \ \ \ \ \ \ \ \ \ \ \ \ \ \left. - ~\frac{(r-1) 
\left.s^{\prime}_{i-1}\right.^2}{w^{\prime}_{i-1}} - 
\frac{r(r-2)s^{\prime}_{i-1}t^{\prime}_{i-1}}{(r-1)b^{\prime}_{i-1}} \right)\\
& t^{\prime}_{i} = t^{\prime}_{i-1} + \pix \left( -2t^{\prime}_{i-1} + 
\frac{2(r-1)\left.s^{\prime}_{i-1}\right.^2}{w^{\prime}_{i-1}}-\frac{2r(r-2)\left.t^{\prime}_{i-1}\right.^2}{(r-1)b^{\prime}_{i-1}} 
\right)\\
&w_0^{\prime}=(1-p_0)^{r+1},~  b_0^{\prime}=rp_0(1-p_0)^r,~ q_0^{\prime 
}=(1-p_0)^{2r},\\
&s_0^{\prime}=(r-1)p_0(1-p_0)^{2r-1},~
t_0^{\prime}=(r-1)^2p_0^2(1-p_0)^{2r-2}
\end{split}
\end{equation}

Note that the auxiliary system of recurrence equations
(\ref{modifiedsystem}) can be converted into a system of
differential equations by means of first order approximations.
Setting $\pix = \epsilon$ in the recurrence equation for
$w^{\prime}_i$ obtained above implies $$w^{\prime}_i -
w^{\prime}_{i-1} = - \epsilon r s^{\prime}_{i-1},$$ so that, for
$\epsilon$ small, we are interested in $\hat{w}, \hat{s}$
satisfying the differential equation
$$\frac{d \hat{w}}{dx} = -r \hat{s}.$$
Applying the same argument to the other recurrence formulae in
(\ref{modifiedsystem}), the following system of differential
equations arises. This system will be referred to as \textit{the
system of differential equations associated with $(r,p_0)$}.
\begin{equation}\label{diffequations}
\begin{split}
&\frac{d \hat{w}}{dx} = -r \hat{s}\\
&\frac{d\hat{b}}{dx} = -\hat{b} - r \hat{t} + r
\hat{s}\\
&\frac{d \hat{q}}{dx} = - \frac{(2r-2)\hat{q}\hat{s}}{\hat{w}}\\
&\frac{d\hat{s}}{dx} = -\hat{s} +
\frac{(r-1)\hat{q}\hat{s}}{\hat{w}} -
\frac{(r-1)\hat{s}^2}{\hat{w}} -
\frac{r(r-2)\hat{s}\hat{t}}{(r-1)\hat{b}}\\
&\frac{d\hat{t}}{dx} = -2\hat{t} + \frac{2(r-1)\hat{s}^2}{\hat{w}}
- \frac{2r(r-2)\hat{t}^2}{(r-1)\hat{b}}\\
&\hat{w}(0)=(1-p_0)^{r+1},~ \hat{b}(0)=rp_0(1-p_0)^r,~ 
\hat{q}(0)=(1-p_0)^{2r},\\
&\hat{s}(0)=(r-1)p_0(1-p_0)^{2r-1},~
\hat{t}(0)=(r-1)^2p_0^2(1-p_0)^{2r-2}.
\end{split}
\end{equation}

Given $p_0 \in (0,1)$, $T>0$ and $\gamma>0$, where  $\gamma <
\min\{w_0,b_0,q_0,s_0,t_0\}$, this system of differential
equations has a solution in the domain $\Omega (\gamma, T) =
\{(x,\hat{w},\hat{b},\hat{q},\hat{s},\hat{t}) \in (-\gamma,T)
\times (\gamma,1)^2 \times (\gamma,1)^3 \}$ which may be uniquely
extended arbitrarily close to the boundary of the domain, by a
standard result in the theory of first order differential
equations (see Hurewicz \cite{hurewicz}, Chapter 2, Theorem 11).

As expected, there is a connection between the original system of
recurrence equations and the system of differential equations
(\ref{diffequations}). This connection is summarised in the lemma
below and follows from the solutions to the original system being
well-approximated by the solutions of the modified system
(\ref{modifiedsystem}), as well as from the relation between the
solutions of (\ref{modifiedsystem}) and of (\ref{diffequations})
given by Euler's method. The proof is routine so is omitted.

\begin{lemma}\label{EulerInt}
Let $r \geq 3$ be an integer and $p_0 \in (0,1)$. Let $k_0>0$ such
that the system of differential equations (\ref{diffequations})
with the initial conditions defined by $p_0$ has positive
solutions in $\Omega$ defined at $x=k_0$. Then, given $\xi>0$,

\begin{enumerate}
\item[(i)] there exists $\epsilon_0>0$ satisfying the following
property. If $0< \epsilon \leq \epsilon_0$ and the system of
recurrence equations (\ref{modifiedsystem}) is solved with $\pix =
\epsilon$, then $|w_i-\hat{w}(\epsilon i)| < \xi$,
$|b_i-\hat{b}(\epsilon i)| < \xi$, $|q_i-\hat{q}(\epsilon i)| <
\xi$, $|s_i-\hat{s}(\epsilon i)| < \xi$ and $|t_i-\hat{t}(\epsilon
i)| < \xi$, for $i = 0, 1, \ldots, \left\lceil k_0/\epsilon
\right\rceil$.

\item[(ii)] there exists $\epsilon_1>0$ such that for $0 \leq
\epsilon \leq \epsilon_1$,
\begin{equation*}
\left|\int_{0}^{k_0} \hat{b}(x)~dx-\sum_{i=0}^{\left\lceil
k_0/\epsilon \right\rceil-1} \epsilon
b_i\right|<\xi,\textrm{ for every $0<\epsilon\leq
\epsilon_1$}.
\end{equation*}
\end{enumerate}

\end{lemma}

Using this lemma, we can now determine additional properties of
the solutions to (\ref{diffequations}).
\begin{lemma}\label{DiffEq}
Given $p_0 \in (0,1)$, the system of differential equations
(\ref{diffequations}) has unique solutions $\hat{w}(x)$, $\hat{b}(x)$,
$\hat{q}(x)$, $\hat{s}(x)$ and $\hat{t}(x)$ defined over the entire
nonnegative real line satisfying the following properties:
\begin{enumerate}
\item[(i)]  $\hat{w}(x)$, $\hat{b}(x)$, $\hat{q}(x)$, $\hat{s}(x)$ and 
$\hat{t}(x)$ are
positive,
\item[(ii)] $\displaystyle{\int_{0}^{\infty} \hat{b}(x)~dx}$ 
converges.
\end{enumerate}
\end{lemma}
\proof
As mentioned before, a standard result in the theory of first
order differential equations ensures that, for $p_0 \in (0,1)$,
$T>0$ and $\gamma < \min\{w_0,b_0,q_0,s_0,t_0\}$, the system of
differential equations has a solution in the domain $\Omega
(\gamma,T) = \{(x,\hat{w},\hat{b},\hat{q},\hat{s},\hat{t}) \in
(-\gamma,T) \times (\gamma,1)^2 \times (\gamma,1)^3 \}$ which may
be uniquely extended arbitrarily close to the boundary of the
domain.

Given $p_0 \in (0,1)$ and $T>0$, we show that there exists
$\gamma=\gamma(T)>0$ such that this system of differential
equations in the domain $\Omega (\gamma,T)$ has a unique solution
defined for $x$ arbitrarily close to $x=T$. This implies that the
solutions are defined over the nonnegative real line.

Suppose on the contrary that, for some $T>0$, no $\gamma(T)$ with
the above property exists. Let $x_0$ denote the infimum of such $T$. 
Let $(x^{\prime}, w^{\prime},
b^{\prime}, q^{\prime}, s^{\prime}, t^{\prime})$ be any point in the
interior of a region $\Omega(\gamma_0,x_0)$ such that a solution to
the system of differential equations exists for $0 \leq x \leq
x^{\prime}$ and $\hat{w}(x^{\prime})=w^{\prime}$,
$\hat{b}(x^{\prime})=b^{\prime}$,
$\hat{q}(x^{\prime})=q^{\prime}$, $\hat{s}(x^{\prime})=s^{\prime}$
and $\hat{t}(x^{\prime})=t^{\prime}$, where $\gamma_0>0$. By
Lemma~\ref{EulerInt}, given $\xi>0$, there exists $\epsilon_0$
such that for $0 < \epsilon \leq \epsilon_0$ and $ 0 \leq i \leq
\lceil x^\prime/\epsilon \rceil$,
$$|w_i-\hat{w}(\epsilon i)| < \xi, |b_i-\hat{b}(\epsilon i)| < \xi, 
|q_i-\hat{q}(\epsilon i)| < \xi, |s_i-\hat{s}(\epsilon i)| < 
\xi,|t_i-\hat{t}(\epsilon i)| < \xi.$$
Since the quantities $w_i$, $b_i$, $q_i$, $s_i$ and $t_i$
represent probabilities of specific events after $i$ steps of a
randomised algorithm, we conclude that
$$q_i + s_i \leq w_i,~ s_i + t_i \leq b_i.$$ Using this and the
fact that $w^{\prime}, b^{\prime}, q^{\prime}, s^{\prime},
t^{\prime}>\gamma_0$ (which is independent of $\xi$), we have
\begin{equation}\label{ineq}
\max\{\hat{q}(x^{\prime}),\hat{s}(x^{\prime})\}<\hat{w}(x^{\prime}),
\quad 
\max\{\hat{s}(x^{\prime}),\hat{t}(x^{\prime})\}<\hat{b}(x^{\prime}).
\end{equation}

Let $m_0$ be a positive integer such that $1/m_0 <
\min\{w_0,b_0,q_0,s_0,t_0\}$. The definition of $x_0$ implies that one 
of the functions $\hat{w}, \hat{b}, \hat{q},
\hat{s}, \hat{t}$ must get arbitrarily close to 0 in the
neighbourhood of a point $x_0$, $0<x_0<T$. By (\ref{ineq}), it
must be one of $\hat{q}$, $\hat{s}$ or $\hat{t}$. (Note that this 
argument also applies in the case that $x=0$.)

Suppose this is the case for $\hat{q}$. Let $x^{\prime}$ be such that the
system of differential equations have a positive solution in
$[0,x^{\prime}]$.  
Recall that
$$\frac{d \hat{q}}{dx} = - \frac{(2r-2)\hat{q}\hat{s}}{\hat{w}} = 
\hat{q} \left(\frac{-(2r-2)\hat{s}}{\hat{w}} \right),$$
and, by equation (\ref{ineq}), $-(2r-2)\hat{s}(x)/\hat{w}(x) \geq
-(2r-2)$ for $0 \leq x \leq x^{\prime}$. Now, if $f$ is the solution for
$$\frac{df}{dx}=-(2r-2)f,~f(0)=\hat{q}(0),$$ we must have $\hat{q}(x) 
\geq f(x)$
for every $x$ in the interval $[0,x^{\prime}]$.  However,
$\displaystyle{f(x)=f(0) e^{-(2r-2)x}}$ is a
strictly positive function in this interval bounded below by the
constant $\displaystyle{f(0) e^{-(2r-2)x^{\prime}}}$.
So $\hat{q}(x)$ cannot approach 0 at $x^{\prime}$. Similar arguments yield
contradictions for the cases when $\hat{s}(x)$ or $\hat{t}(x)$
approach 0 in the neighbourhood of the point $x_0$, since
\begin{equation*}
\begin{split}
\frac{d\hat{s}}{dx} &= \left(-1 +
\frac{(r-1)\hat{q}}{\hat{w}} -
\frac{(r-1)\hat{s}}{\hat{w}} -
\frac{r(r-2)\hat{t}}{(r-1)\hat{b}}\right)\hat{s}\\
&\geq \left(-1-(r-1)-\frac{r(r-2)}{r-1}\right) \hat{s},
\end{split}
\end{equation*}
and
\begin{equation*}
\begin{split}
\frac{d\hat{t}}{dx} &= -2\hat{t} + \frac{2(r-1)\hat{s}^2}{\hat{w}}
- \frac{2r(r-2)\hat{t}^2}{(r-1)\hat{b}}\\
& \geq \left(-2-\frac{2r(r-2)}{r-1}\right) \hat{t}.
\end{split}
\end{equation*}

Thus, the solutions to the system of differential equations are
indeed defined over the entire nonnegative real line. Furthermore,
the previous argument ensures that they are positive, concluding
the proof of part (i).

For part (ii), note that the differential equations for $\hat{w}$
and $\hat{b}$ in (\ref{diffequations}) imply
$$\frac{d(\hat{w}+\hat{b})}{dx}=-\hat{b}-r\hat{t},$$ so
$$\hat{b}(x) = -\frac{d(\hat{w}+\hat{b})}{dx}(x)-r\hat{t}(x) \leq
-\frac{d(\hat{w}+\hat{b})}{dx}(x),~\forall x.$$ As a consequence,
for every $T>0$,
$$\int_{0}^T \hat{b}(x)~dx \leq \int_0^{T} 
-\frac{d(\hat{w}+\hat{b})}{dx}(x)~dx = \hat{w}(0) +\hat{b}(0) - \hat{w}(T) - \hat{b}(T) \leq  
\hat{w}(0) +\hat{b}(0).$$
This proves part (ii). \qed

\section{Proof of Theorem~\ref{maintheorem}}\label{finalsection}

As mentioned in the introduction, we wish to obtain a lower bound
on the cardinality of a largest vertex subset that induces a
forest in an $r$-regular graph $G$ not containing short cycles.
Recall our definition of $\tau(G)$, given by
\begin{equation*}
\tau(G)=\max~\{|V(F)| ~:~ F \textrm{ is an induced forest in }G\}.
\end{equation*}

Let $G$ be an $r$-regular graph on $n$ vertices with girth $g$ and
consider the set $P$ of purple vertices at the end of step 2 when
Algorithm~\ref{alg1} is applied to $G$ with $N<g/2-2$. It
is clear that the induced graph $G[P]$ contains a cycle only if
some vertex $v$ with at least two purple neighbours has been added
to $P$. By the description of our algorithm, this cannot happen
unless $v$ was selected in the same step as one of its neighbours.
It follows that, if $\bar{P}$ is the set obtained from $P$ by
deleting any pairs of adjacent vertices added to $P$ in the same
step, the induced subgraph $G[\bar{P}]$ is acyclic.

Now, given a vertex in $R_i$, the probability that none of its
neighbours is also selected is at least $(1-\pix )^r$, since a
vertex has at most $r$ neighbours that could be added to $R_i$.
Therefore, the expected number of vertices added to $P$ at time
$i$ that are not removed is at least $\pix (1-\pix )^rb_{i-1}n$
and
\begin{equation}\label{sizeP}
{\bf E}|\bar{P}| \geq p_0(1-p_0)^r n + \sum_{i=1}^N \pix
(1-\pix )^rb_{i-1}n
\end{equation}

Part of the set $W$ of white vertices produced at the end of the
algorithm will also be added to the forest. By definition, these
vertices have no purple neighbours, so that no cycle containing
purple vertices is created by adding white vertices to $\bar{P}$.
Thus $G[\bar{P} \cup \bar{W}]$ is still acyclic, where $\bar{W}$
denotes the set of vertices in acyclic components of $G[W]$.

Now, since $G$ has girth $g$, no cycles appear if we add white
vertices lying in components of $G[W]$ of size at most $g-1$.
Therefore, a lower bound on the size of $\bar{W}$ can be obtained
by estimating the number of vertices in small components of
$G[W]$. This will be done through a branching process argument.

To define the branching process, start with a white vertex $v_0$
and set the random variable $Y_0=\{v_0\}$. In general, $Y_i$
denotes the set of white vertices already exposed, but whose
neighbours have not been considered yet. Define $U_0=V(G)-\{v_0\}$
and let $U_i$ be the random variable accounting for the set of
vertices which have not been exposed by the branching process up
to step $i$. After step $i$, either $|Y_i|=0$, in which case the
process has died out, or $|Y_i|>0$, in which case we choose a
white vertex $v_i$ in $Y_i$, expose its white neighbours
$\displaystyle{N_W(v_i)\subseteq U_i}$ and define
$\displaystyle{Y_{i+1}=Y_i \cup N_W(v_i) - \{v_i\}}$, $U_{i+1} =
U_i \setminus N(v_i)$. We are interested in estimating the
probability that $|Y_{g-1}|>0$, i.e., that the branching process
has not died out after $g-1$ steps.

\begin{prop}\label{branch}
Let $\delta>0$, fix an integer $r\geq 3$ and suppose the existence
of $p_0>0$ such that the solutions to the system of differential
equations associated with $(r,p_0)$ satisfy
$\displaystyle{\lim_{x \rightarrow \infty}
\frac{(r-1)\hat{q}(x)}{\hat{w}(x)} < 1}$. Then, there exist $g>0$,
$0<N<g/2-1$ and $0<p<1$ such that, if Algorithm~\ref{alg1} is
applied to an $r$-regular graph $G$ with girth at least $g$ for
$N$ steps with probabilities $(p_0,p)$, then
$$\pr(|Y_{g/2-1}|>0) < \delta.$$
\end{prop}

\proof Let $Z_i$ denote the random variable counting the number of
neighbours of $v_i$ in $U_i$. Note that $Z_0$ has binomial
distribution ${\bf Bin}(r,q_N/w_N)$, since 
Corollary~\ref{lem3} and
Lemma~\ref{lem2} ensure that, conditional upon $v_0$ being white,
the events associated with each of its neighbours being white are
mutually independent and have probability $q_N/w_N$. Furthermore,
$Z_i$ has distribution ${\bf Bin}(r-1,q_N/w_N)$ for every $i \geq 1$,
since the condition $0<i<g-1$ implies $|N(v_i) \cap U_i|=r-1$, and
Corollary~\ref{lem3} and Lemma~\ref{lem2} are applicable in the
same way.

Let $k_0>0$ such that the solution to the system of differential
equations $(\ref{diffequations})$ satisfies
$(r-1)\hat{q}(k_0)/\hat{w}(k_0) < 1$.

Let $\xi >0$ be such that
$\displaystyle{\frac{(r-1)(\hat{q}(k_0)+\xi)}{\hat{w}(k_0)-\xi} <
1}$. Fix $\epsilon_0$ as in Lemma~\ref{EulerInt}, part $(i)$, and
let $\epsilon<\epsilon_0$ such that $N=k_0/\epsilon$ is an
integer. Now, apply Algorithm~\ref{alg1} for $N$ steps with the
given $p_0$ and $\pix = \epsilon$, for all $i \geq 1$, to a graph
$G$ with girth $g \geq 2N+3$. Then,
\begin{equation*}
\frac{(r-1)q_N}{w_N} \leq \frac{(r-1)(\hat{q}(k_0) + 
\xi)}{\hat{w}(k_0)-\xi} < 1
\end{equation*}
So, we have $(r-1)q_N/w_N<1$, and a branching process argument as
in \cite{alon_spencer} shows that, by choosing $g$ sufficiently large, 
$\pr(|Y_{g/2-1}|>0)< \delta$, as
required. \qed

By the above proposition, given $\delta>0$ and $p_0>0$ such that
the solutions to the system of differential equations associated
with $(r,p_0)$ satisfy $\displaystyle{\lim_{x \rightarrow
\infty} \frac{(r-1)\hat{q}(x)}{\hat{w}(x)} < 1}$, we may fix $g$,
$N$ and $p$ so as to have the property $\pr(Y_{g/2-1}>0) <
\delta$, i.e., $\pr(Y_{g/2-1}=0) \geq 1 - \delta$. It follows that
for such $g$ the expected number of white vertices in acyclic
components of $G[W]$ is bounded below by
\begin{equation}\label{sizeW}
(1 - \delta) w_N n.
\end{equation}

We are now ready to prove Theorem~\ref{maintheorem}.
\medskip

\noindent \textbf{Proof of Theorem~\ref{maintheorem}\ } Fix $r \in
\mathbb{N}$ and $\delta>0$. We show that, given $p_0 \in (0,1)$,
the inequality  $\displaystyle{\tau(G) \geq
\left(\xi(p_0)-\delta\right)n}$ holds, where
\begin{equation*}
\begin{matrix}
\\
\xi(p_0) =\\
\\
\end{matrix}
\left\{
\begin{matrix}
\displaystyle{p_0(1-p_0)^r + \int_{0}^{\infty} \hat{b}(x)~dx
 + \lim_{x \rightarrow \infty} \hat{w}(x)}, &  \textrm{if } 
\displaystyle{\lim_{x \rightarrow \infty} \frac{(r-1)\hat{q}(x)}{\hat{w}(x)}<1}  
\\
\displaystyle{p_0(1-p_0)^r + \int_{0}^{\infty} \hat{b}(x)~dx},&
\textrm{otherwise.}\ \ \ \ \ \ \ \ \ \ \ \ \ \ \ \
\end{matrix}
\right.  
\end{equation*}
Here, $\hat{w}$, $\hat{b}$ and $\hat{q}$ are solutions to the
system of differential equations associated with $(r,p_0)$. By
Lemma~\ref{DiffEq}, this system has positive solutions
$\hat{w},\hat{b},\hat{q},\hat{s},\hat{t}$ defined over the
nonnegative real line such that $\displaystyle{\int_{0}^{\infty}
\hat{b}(x)~dx}$ converges.

Let $k_0>0$ be such that, for every $k>k_0$,
\begin{equation}\label{fix k}
\left|\int_{0}^{\infty} \hat{b}(x)~dx-\int_{0}^{k}
\hat{b}(x)~dx\right|<\frac{\delta}{6}
\end{equation}
Using Lemma~\ref{EulerInt}, fix $\epsilon_0>0$ such that
$$|w_i-\hat{w}(\epsilon i)| < \frac{\delta}{6},~i = 0, 1, \ldots, 
\left\lceil \frac{k_0}{\epsilon} \right\rceil$$
and fix $\epsilon_1>0$ satisfying
\begin{equation*}
\left|\int_{0}^{k_0} \hat{b}(x)~dx-\sum_{i=0}^{\left\lceil
\frac{k_0}{\epsilon} \right\rceil-1} \epsilon
b_i\right|<\frac{\delta}{6},\textrm{ for every $0<\epsilon\leq
\epsilon_1$}.
\end{equation*}

Let $\epsilon=\min \{\epsilon_0, \epsilon_1 , 1-(1-\delta/6)^{1/r}
\}$ and $N=\left\lceil k_0/\epsilon \right\rceil$. Fix $g = 2N +
3$. Then, given an $r$-regular graph  $G$ with girth larger than
or equal to $g$, we apply Algorithm~\ref{alg1} for $N$ steps with
probabilities $(p_0,p=\epsilon)$. The first moment principle leads
to a lower bound for $\tau(G)$.  As a matter of fact, our lower
bound (\ref{sizeP}) on the cardinality of $\bar{P}$ implies
\begin{equation}\label{sizeP2}
\begin{split}
{\bf E}|\bar{P}| &\geq np_0(1-p_0)^r +
n(1-\epsilon)^r\left(\sum_{i=1}^N\epsilon
b_{i-1}\right)\\
& \geq np_0(1-p_0)^r + n(1-\epsilon)^r\left(\int_{0}^{k_0}
\hat{b}(x)~dx - \frac{\delta}{6}\right)\\
& \geq np_0(1-p_0)^r + n\left(1 -
\frac{\delta}{6}\right)\left(\int_{0}^{\infty}
\hat{b}(x)~dx -\frac{2\delta}{6}\right)\\
& \geq n\left(p_0(1-p_0)^r + \int_{0}^{\infty} \hat{b}(x)~dx
\right) - \frac{\delta~n}{2},
\end{split}
\end{equation}

If, in addition, the solutions to the system of differential
equations associated with $(r,p_0)$ satisfy
$\displaystyle{\lim_{x \rightarrow \infty}
\frac{(r-1)\hat{q}(x)}{\hat{w}(x)} < 1}$, Proposition~\ref{branch}
establishes a lower bound (\ref{sizeW}) on the cardinality of the
set $\bar{W}$ of white vertices that can be added to the forest.
Clearly, $k_0$ in (\ref{fix k}) may be chosen so that, for every
$k>k_0$, we also have
$$\left|\lim_{x->\infty}\hat{w}(x)- 
\hat{w}(k)\right|<\frac{\delta}{6}$$
and
$$\frac{(r-1)\hat{q}(k)}{\hat{w}(k)}<1.$$
The girth $g$ can also be taken larger, if necessary, to ensure
that the size of $\bar{W}$ is bounded below by
$\displaystyle{(1-\delta/6) w_N n}$.

Thus,
\begin{equation}\label{sizeW2}
\begin{split}
{\bf E}|\bar{W}| & \geq  n\left(1 - \frac{\delta}{6}\right)w_N  \geq 
n\left(1 - \frac{\delta}{6}\right)^2\hat{w}(\epsilon N)\\
& \geq n\left(1 - \frac{\delta}{6}\right)^3 \lim_{x \rightarrow
\infty} \hat{w}(x) \geq n\lim_{x \rightarrow \infty} \hat{w}(x) -
\frac{\delta~n}{2},
\end{split}
\end{equation}

Now, given that $\tau(G) \geq {\bf E}|\bar{P} \cup \bar{W}|$ and
using equations (\ref{sizeP2}) and (\ref{sizeW2}), we conclude
that
$$\tau(G) \geq \left(\xi(p_0)-\delta\right)n,$$ as claimed.
Numerical calculations of these quantities lead us to the bounds
in Table~1. We note that, for every value of $r$ tested, we were
able to choose a constant $p_0$ such that the numerical solutions
to the system of differential equations associated with $(r,p_0)$
satisfy $\displaystyle{\lim_{x \rightarrow \infty}
\frac{(r-1)\hat{q}(x)}{\hat{w}(x)} < 1}$. \qed

\end{document}